\documentclass[noamsfonts]{amsproc}
\usepackage{aspmproc, amssymb}

\title[On the cusp form motives in genus~$1$ and level~$1$]{On the cusp 
form motives in genus~$1$ and level~$1$}

\author[C. Consani and C. Faber]{Caterina Consani and Carel Faber}

\address{
{\rm Caterina Consani}\\
Department of Mathematics \hfill Department of Mathematics\\
The Johns Hopkins University \hfill University of Toronto\\
3400 North Charles Street \hfill Toronto, Ontario\\
Baltimore, MD 21218, U.S.A. \hfill Canada M5S 3G3\\
kc@math.jhu.edu \hfill kc@math.toronto.edu\\
\\
{\rm Carel Faber}\\
Institutionen f\"or Matematik \hfill Department of Mathematics\\
Kungliga Tekniska H\"ogskolan \hfill The Johns Hopkins University\\
100 44 Stockholm \hfill 3400 North Charles Street\\
Sweden \hfill Baltimore, MD 21218, U.S.A.\\
faber@kth.se \hfill faber@math.jhu.edu}

\rcvdate{April 20, 2005}
\rvsdate{December 23, 2005}

\thanks{This paper owes a great debt to the work of
Ezra Getzler. We thank Mathematisches Forschungsinstitut Oberwolfach, where
we began this work during a ``Research in Pairs'' stay in May 2004.
This research was supported in part by NSERC grant 72024520 and by grants from
the Swedish Research Council and the G\"oran Gustafsson Foundation.}

\newtheorem{pr}{Proposition}
\newtheorem{lm}{Lemma}
\newtheorem{tm}{Theorem}

\newcommand{\D}{{\mathcal{D}}}
\newcommand{\E}{{\mathcal{E}}}
\newcommand{\SSS}{\mathsf{S}}

\newcommand{\Q}{\mathbb{Q}}

\newcommand{\Z}{\mathbb{Z}}

\newcommand{\lan}{\langle}
\newcommand{\ran}{\rangle}

\newcommand{\bpf}{\noindent {\em Proof.} }
\newcommand{\epf}{\qed \vspace{+10pt}}

\newcommand{\mbar}[1]{{\overline{M}}_{#1}}
\newcommand{\m}[1]{M_{#1}}
\newcommand{\bm}[1]{{\partial M}_{#1}}
\newcommand{\Sn}{\Sigma_n}
\newcommand{\sn}{\Sigma_{n-1}}
\newcommand{\bA}[1]{{\bf a}_{#1}}
\newcommand{\bB}[1]{{\bf b}_{#1}}
\newcommand{\Alt}{{\rm{Alt}}}
\newcommand{\p}{\partial}
\newcommand{\ch}{{\rm{ch}}}
\newcommand{\Ch}{{\rm{Ch}}}
\newcommand{\sgn}{{\rm{sgn}}}
\def\Lie{{\operatorname{\mathcal{L}\mathit{ie}}}}
\def\({(\!(}
\def\){)\!)}

\newcommand{\De}[1]{\Delta_{#1}}
\newcommand{\Do}[1]{\Delta^{\circ}_{#1}}
\newcommand{\Ind}{{\rm Ind}}
\newcommand{\Sym}{{\rm Sym}}
\newcommand{\pp}{{\rm pr}}

\begin{document}

\begin{abstract}
We prove that the moduli space of stable $n$-pointed curves of genus~$1$
and the projector associated to the alternating representation of the
symmetric group on $n$ letters define (for $n>1$) the Chow motive
corresponding to cusp forms of weight $n+1$ for ${\rm SL}(2,\Z)$.
This provides an alternative (in level~$1$) to the construction of Scholl.
\end{abstract}

\maketitle

\setcounter{section}{0}
\section{Introduction}
In this paper we give an alternative construction of the Chow motives
$S[k]$ corresponding to cusp forms of weight $k$ for ${\rm SL}(2,\Z)$.
The Betti cohomology related to these cusp forms was initially studied
by Eichler and Shimura, after which Deligne constructed the corresponding
$\ell$-adic Galois representations. Using the canonical desingularization
of the fiber products of the compactified universal elliptic curve
constructed by Deligne, Scholl then defined projectors 
such that the realizations of the associated Chow motives are
these parabolic cohomology groups.
The smooth projective varieties
used in this construction are called Kuga-Sato varieties.

Instead of the Kuga-Sato varieties, we use the spaces $\mbar{1,n}$, the
Knudsen-Deligne-Mumford moduli spaces of stable $n$-pointed curves
of genus~$1$. The symmetric group $\Sn$ acts naturally on $\mbar{1,n}$,
by permuting the $n$ marked points. Let $\alpha$ denote its alternating 
character. Our main result is that $\mbar{1,n}$ (for $n>1$) and the projector
$\Pi_{\alpha}$ corresponding to $\alpha$ define the Chow motive $S[n+1]$.
In other words, we have the following result.
\par\medskip
\noindent {\bf Theorem.} {\em 
For $n>1$,}
$$\Pi_{\alpha}(H^*(\mbar{1,n},\Q))=\Pi_{\alpha}(H^n(\mbar{1,n},\Q))
=H^1_!(\m{1,1},\Sym^{n-1}R^1\pi_*\Q).$$
\noindent Here $\pi:\E\to\m{1,1}$ is the universal elliptic curve and
$H^i_!={\rm Im}(H^i_c\to H^i)$ denotes the parabolic cohomology.

The cohomology $H^*(\mbar{g,n})$ of the moduli space of stable $n$-pointed
curves of genus~$g$ has been studied intensively in recent years, in
particular for $n>0$ through the connection with Gromov-Witten theory.
Since $\mbar{g,n}$ is a smooth projective stack over $\Z$, these groups
have arithmetic relevance as well. Getzler has initiated the study of the
cohomology $H^*(\m{g,n})$ of the moduli space of smooth $n$-pointed curves
of genus~$g$ as a representation of $\Sn$. Through the theory of
modular operads, as developed by Getzler and Kapranov, the $\Sn$-equivariant
Euler characteristics of the cohomology of 
the spaces $\mbar{g,n}$ are expressed
in the $\Sn$-equivariant Euler characteristics of the cohomology of 
the spaces $\m{g,n}$. The action of $\Sn$ is 
crucial here. Another central idea
of Getzler is to express the $\Sn$-equivariant Euler characteristic of
$H^*(\m{g,n})$ in terms of the Euler characteristics of the cohomology
of irreducible symplectic local systems on $\m{g}$. Since these local systems
are pulled back from the moduli space $A_g$ of principally polarized
abelian varieties of dimension~$g$, this provides a connection with genus~$g$
Siegel modular forms.

In genus~$1$, this connection is given by Eichler-Shimura theory. In higher
genus, despite the very important work of Faltings and Chai, much less is
known. Van der Geer and the second author have obtained an explicit
conjectural formula for the motivic Euler characteristics of these local systems
in genus~$2$.

Our work is motivated by the desire to understand the motives underlying
Siegel modular forms and the cohomology of the corresponding local systems.
We expect that the results proved in this paper for genus~$1$, when
suitably generalized, will provide a major step towards this goal.

Unbeknownst to us, Manin had suggested in \cite{Ma1},~0.2 and~2.5, that
it would be desirable to replace the Kuga-Sato varieties by
moduli spaces of curves of genus~$1$ with marked points and a level
structure. Cf.~\cite{Ma2},~3.6.2.

In section~2, we determine the alternating part of the $\Sn$-equivariant
Euler characteristic of $\m{1,n}$. Section~3 deals with the
$\Sn$-equivariant cohomology of $\m{0,n}$; some of the results obtained
here may be of independent interest. The theory of modular operads
and the results of section~3 are used in section~4 to determine the
alternating part of the Euler characteristic of $\mbar{1,n}\setminus\m{1,n}$.
In section~5 we combine the results of sections~2 and 4 and prove our main
theorem.

\section{The contribution of the interior}
In this section we determine the contribution of
$\m{1,n}$,
i.e., we determine
$$\lan s_{1^n},e_c^{\Sn}(\m{1,n})\ran,$$
the alternating part of the $\Sn$-equivariant Euler characteristic 
of the compactly supported cohomology
of $\m{1,n}$.
Here, for a partition~$\lambda$ of~$n$, the notation~$s_{\lambda}$ is used
for the Schur function corresponding to the irreducible representation
of~$\Sn$ indexed by~$\lambda$, and~$\lan\,,\,\ran$ stands for the
standard inner product on the ring of symmetric functions, for which
the~$s_{\lambda}$ form an orthonormal basis.
We will usually not make a notational distinction between a (possibly 
virtual) $\Sn$-representation~$V$ and its characteristic~$\ch_n(V)$, 
the symmetric function corresponding to it (\cite{GK},~7.1).
The Euler characteristic is taken in $K_0$ of a convenient category, such as
the category of mixed Hodge structures or of $\ell$-adic Galois representations.

Let $E=(E,0)$ be an elliptic curve. We may think of the points
of $E^{n-1}$ as $n$-tuples
$$(0,x_2,\dots,x_n)$$
(with $x_1=0$)
and by doing so we find a natural action of $\Sn$ on $E^{n-1}$
(combine the effect of a permutation $\sigma$ with a translation
of each coordinate over $-x_{\sigma^{-1}(1)}$).
We are interested in the subspace of $H^\bullet(E^{n-1})
=H^\bullet(E)^{\otimes(n-1)}$ where the induced action of $\Sn$
is via the alternating representation.

Let $\sn\subset\Sn$ be the subgroup permuting the last $n-1$ entries.
\begin{lm} The subspace of $H^\bullet(E)^{\otimes(n-1)}$
where the induced action of $\sn$ is
via the alternating representation is isomorphic to
$$\oplus_{k=0}^{n-1}
\wedge^k H^{{\rm even}}(E)\otimes\Sym^{n-1-k}H^1(E).$$
\end{lm}

\bpf
The subspace of $V=H^\bullet(E)^{\otimes(n-1)}$ where $\sn$ acts
alternatingly is generated by sums
$$\sum_{\sigma\in\sn}(-1)^{\sgn(\sigma)}\sigma^*(v)$$
with $v\in V$.
Clearly, we may restrict ourselves to pure tensors $v$ such that the first
$k$ factors are in $H^{{\rm even}}(E)$ and the remaining $n-1-k$ factors
are in $H^1(E)$, for some $k$. Fix $k$. It suffices now to consider the
action of $\Sigma_k\times\Sigma_{n-1-k}$ on such $v$. This leads to the
claimed isomorphism.\epf

Note that only the terms with $k\le2$ in the direct sum above are
nonzero.
Thus it is concentrated in degrees $n-2$, $n-1$, and $n$.
\begin{pr} The subspace of $H^\bullet(E)^{\otimes(n-1)}$ where
the induced action of $\Sn$ is via the alternating representation
is $\Sym^{n-1}H^1(E)$.
\end{pr}

\bpf Let $\tau\in\Sn$ be the transposition $(12)$. We need to show
that $\tau^*(\gamma)=-\gamma$ for all
$\gamma\in\Sym^{n-1}H^1(E)$, but that none of the
$\sn$-alternating vectors 
coming from 
elements of $H^{{\rm
even}}(E)\otimes\Sym^{n-2}H^1(E)$ and $\wedge^2 H^{{\rm
even}}(E)\otimes\Sym^{n-3}H^1(E)$ have this property.

As an example, consider the case $n=2$. Note that
$\tau(0,x_2)=(x_2,0)=(0,-x_2)$. Thus $\tau=-1_E$ and the
$(-1)$-eigenspace of $\tau^*$ on $H^\bullet(E)$ is $H^1(E)$.

In the general case,
$$\tau(0,x_2,x_3,\dots,x_n)=(0,-x_2,x_3-x_2,\dots,x_n-x_2).$$
Denote by $\pp_i:E^{n-1}\to E$ the projection onto the $i$th
factor (with $2\le i\le n$) and by $\tau_i$ the composition
$\pp_i\circ\tau$. Then
$$\tau^*(\gamma_2\otimes\dots\otimes\gamma_n)=
\tau_2^*(\gamma_2)\cdot\ldots\cdot\tau_n^*(\gamma_n).$$ For
$k\ge3$ we have $\tau_k=m\circ((-\pp_2)\times\pp_k)$, where
$m:E\times E\to E$ denotes the group law. Observe now that
$$\tau_k^*(\zeta)=-\pp_2^*(\zeta)+\pp_k^*(\zeta)$$
for $\zeta\in H^1(E)$ and $k\ge3$.

Denote by $p_i:E^{n-1}\to E^{n-2}$ the projection forgetting the
$i$th factor ($2\le i\le n$). Let $\gamma\in\Sym^{n-2}H^1(E)$ and
denote by
$$\Gamma=\sum_{i=2}^n (-1)^i p_i^*\gamma$$
the $\sn$-alternating vector corresponding to $1\otimes\gamma$.
Let $I$ be the ideal $\pp_2^*(H^1(E)\oplus H^2(E))$. Note that
$\Gamma\equiv p_2^*\gamma=1\otimes\gamma\mod I$. But
$\tau^*(1\otimes\gamma)\equiv 1\otimes\gamma\mod I$ by the above.
Thus $\tau^*\Gamma=-\Gamma$ implies $\gamma=0$.

This shows that the $\sn$-alternating vectors corresponding to
elements of $H^0(E)\otimes\Sym^{n-2}H^1(E)$ are not
$\Sn$-alternating. We conclude that the alternating representation
of $\Sn$ does not occur in degree $n-2$. By duality, it does not
occur in degree $n$ either.

Denote by $p_{ij}:E^{n-1}\to E^{n-3}$ the projection forgetting
the $i$th and $j$th factors ($2\le i<j\le n$). Let
$\gamma\in\Sym^{n-3}H^1(E)$ and denote by
$$\Xi=\sum_{i=2}^{n-1}\sum_{j=i+1}^n (-1)^{i+j}
p_{ij}^*\gamma\cdot(\pp_i^*p-\pp_j^*p)$$ the $\sn$-alternating
vector corresponding to $(1\wedge p)\otimes\gamma$ (here $p$ is
the class of a point). Then
$$(-1)^{n+1}p_{n*}\Xi=p_{n*}\big(\sum_{i=2}^{n-1}(-1)^i
p_{in}^*\gamma\cdot\pp_n^*p\big) =\sum_{i=2}^{n-1} (-1)^i
p_i^*\gamma,$$ the $\Sigma_{n-2}$-alternating vector in
$H^\bullet(E)^{\otimes(n-2)}$ corresponding to $1\otimes\gamma$.
Using that $p_n\circ\tau=\tau\circ p_n$, one shows that
$\tau^*\Xi=-\Xi$ implies $\gamma=0$. Thus the alternating
representation of $\Sn$ can occur only in $\Sym^{n-1}H^1(E)$.

To conclude, 
we show that these vectors are indeed $\Sn$-alternating.
Choose $\alpha$ and $\beta$ in $H^1(E)$ with $\alpha\cdot\beta=p$. Fix
$k$ and $l$ with sum $n-1$ and let
$\gamma=\gamma_2\otimes\dots\otimes\gamma_n$
with $k$ of the factors equal to $\alpha$ and the remaining $l$
equal to $\beta$.
If $\gamma_2=\alpha$, then $\gamma+\tau^*\gamma$ is a sum of $l$ terms;
each term arises from $\gamma$ by replacing $\gamma_2$ by $p$ and one
of the $\beta$'s by $1$. If $\gamma_2=\beta$, then $\gamma+\tau^*\gamma$
is a sum of $k$ terms; each term arises from $\gamma$ by
replacing $\gamma_2$ by $-p$ and one of the $\alpha$'s by $1$. It is
now easy to see that the symmetric tensor $\Gamma$ that is the sum of
all $\gamma$ satisfies $\tau^*\Gamma=-\Gamma$. This finishes the proof.\epf

We may think of the fiber of $\m{1,n}$ over $[E]$ as the open
subset $D^\circ_n$ of $E^{n-1}$ where the $n$ points $0,x_2,\dots,x_n$
are mutually distinct, i.e., the complement of the $n-1$ zero
sections $x_i=0$ (with $2\le i\le n$) and the diagonals $x_i=x_j$
(with $2\le i<j\le n$). Clearly this open subset is $\Sn$-invariant.
\begin{lm}
The subspace of $H^\bullet_c(D^\circ_n)$
where the induced action of $\Sn$ is
via the alternating representation is canonically isomorphic to
the corresponding subspace of $H^\bullet(E^{n-1})$, thus to
$\Sym^{n-1}H^1(E)$.
\end{lm}

\bpf Write $D_k$ for the closed subset of $E^{n-1}$ where $\{0,x_2,\dots,x_n\}$
has cardinality at most $k$ and $D_k^\circ=D_k\setminus D_{k-1}$ for its
open subset where $\{0,x_2,\dots,x_n\}$ has cardinality $k$. The subsets
$D_k$ and $D_k^\circ$ are $\Sn$-invariant. By induction on $k$ we show that
$H^\bullet_c(D_k)$ does not contain a copy of the alternating representation
for $k\le n-1$. Note that $D_1$ is a point. We may assume $n>2$. We have
exact sequences
$$H^{i-1}_c(D_{k-1})\to H^i_c(D_k^\circ) \to H^i_c(D_k) \to H^i_c(D_{k-1})$$
of $\Sn$-representations. By induction, the outer terms do not contain
alternating representations. Consider $D_k^\circ$ for $k\le n-1$. For every
connected component, there exists 
a transposition
in $\Sn$ acting on it as the
identity. This shows that $H^\bullet_c(D_k^\circ)$ does not contain an
alternating representation, and the same holds for $H^\bullet_c(D_k)$.
The exact sequence above, with $k=n$, now gives the result.\epf

For a variety $X$ with $\Sn$-action, denote $\lan s_{1^n},e_c^{\Sn}(X)\ran$
by $A_c(X)$. Clearly, we have
$$A_c(E^{n-1})=A_c(D_n^\circ)=(-1)^{n-1}\Sym^{n-1}H^1(E).$$
Let $\pi:\E\to S$ be a relative elliptic curve. We may consider the
$\Sn$-action on the relative spaces $\E^{n-1}/S$ and $\D_n^\circ/S$
and obtain
$$A_c(\E^{n-1}/S)=A_c(\D_n^\circ/S)=(-1)^{n-1}\Sym^{n-1}R^1\pi_*\Q$$
and similarly with $\Q_{\ell}$-coefficients. The Leray spectral sequence
gives then immediately
$$A_c(\E^{n-1})=A_c(\D_n^\circ)=(-1)^{n-1}e_c(S,\Sym^{n-1}R^1\pi_*\Q).$$
Applying this to the universal elliptic curve, we obtain in particular
$$A_c(\m{1,n})=(-1)^{n-1}e_c(\m{1,1},\Sym^{n-1}R^1\pi_*\Q).$$

Let $n>1$. Then $H^i_c(\m{1,1},\Sym^{n-1}R^1\pi_*\Q)=0$ when $i\neq1$
or $n$ even. For $n$ odd,
$$H^1_c(\m{1,1},\Sym^{n-1}R^1\pi_*\Q)=S[n+1]+1,$$
cf.~\cite{Ge4}, Thm.~5.3 and below. Here we have written $1$ for the
trivial Hodge structure $\Q$ (or the corresponding $\ell$-adic Galois
representation) and $S[n+1]$ for Getzler's $\SSS_{n+1}$; this is an
equality in the Grothendieck group of our category.
We have proved the following result.
\begin{tm}
\label{interior}
The alternating part of the $\Sn$-equivariant Euler characteristic
of the compactly supported cohomology of~$\m{1,n}$ is given by the
following formula:
$$A_c(\m{1,n})=
\left\{
\begin{array}{rl}
-S[n+1]-1,&\quad n>1\,\,{\rm odd;}\\
0,&\quad n\,\,{\rm even.}\\
\end{array}
\right.
$$
Here $S[n+1]=H^1_!(\m{1,1},\Sym^{n-1}R^1\pi_*\Q)$, the parabolic
cohomology of the local system~$\Sym^{n-1}R^1\pi_*\Q$, is the part
of the cohomology of~$\m{1,n}$ corresponding to cusp forms of weight~$n+1$.
\end{tm}
Of course $A_c(\m{1,1})=L$, the Hodge structure $\Q(-1)$. If we
formally define $S[2]=-L-1$, then the formula above holds for
$n=1$ as well.

\section{Cohomology of genus~$0$ moduli spaces and representations
of symmetric groups}

In this section we study the cohomology groups $H^i(\m{0,n})$
as representations of the symmetric group $\Sn$. One of our main tools
is the following. Let $X$ be an algebraic variety, let $Y\subset X$
be a closed subvariety, and let $U=X\setminus Y$ denote the complement. Then
the long exact sequence of compactly supported cohomology
$$\dots\to H^k_c(U)\to H^k_c(X)\to H^k_c(Y)\to H^{k+1}_c(U)\to\dots$$
is a sequence of mixed Hodge structures. See~\cite{DK},~p.~282.

\begin{lm} {\rm (Getzler)}
\label{wt2i}
The mixed Hodge structure on $H^i(\m{0,n})$ is pure of weight $2i$.
\end{lm}

\bpf This is Lemma~3.12 in \cite{Ge1}. We wish to give a different proof
here. The case $n=3$ is trivial. For $n=4$, we use the sequence above,
with $X={\bf P}^1$, $Y=\{0,1,\infty\}$, and $U=\m{0,4}$. The sequence reads
$$0\to H^0_c(\m{0,4})\to H^0_c({\bf P}^1) \to H^0_c(\{0,1,\infty\})
\to H^1_c(\m{0,4})\to 0 \to 0 \to $$
$$\to H^2_c(\m{0,4})\to H^2_c({\bf P}^1) \to 0. $$
Note first that $H^0_c(\m{0,4})=0$. Clearly, $H^1_c(\m{0,4})$ has weight $0$
and $H^2_c(\m{0,4})$ has weight $2$. The statement follows by duality:
$$H^k_c(V)^{\vee}\cong H^{2m-k}(V)(m)$$
as mixed Hodge structures, for $V$ a nonsingular irreducible 
variety of dimension~$m$. 

For $n>4$ we have that $U=\m{0,n}$ is isomorphic to the complement in
$X=\m{0,n-1}\times \m{0,4}$ of the disjoint union
$$Y=\coprod_{i=4}^{n-1} \{x_i=x_n\}\,,$$
where we think of a $k$-pointed curve of genus~$0$ as given by a $k$-tuple
$(0,1,\infty,x_4,\dots,x_k)$ on ${\bf P}^1$. Thus,
$$H^{k-1}_c(Y)\to H^k_c(U)\to H^k_c(X)$$
is an exact sequence of mixed Hodge structures. By dualizing and 
applying a Tate twist,
the same holds for
$$H^i(X)\to H^i(U) \to H^{i-1}(Y)(-1)$$
(with $i=2(n-3)-k$).
By the K\"unneth formula and induction on $n$, the terms on the left and right
have pure Hodge structures of weight $2i$. Hence the same holds for the
term in the middle.
\epf

For $k\ge0$, denote by $\De{k}$ the closed part of $\mbar{0,n}$
corresponding to stable curves with at least $k$ nodes and denote
by $\Do{k}$ the open part $\De{k}\setminus\De{k+1}$ corresponding to
stable curves with exactly $k$ nodes. Put $d=n-3$. Clearly,
$\De{k}\neq\emptyset$ for $0\le k\le d$. In general, $\De{k}$ is
singular, with nonsingular irreducible components, all of
codimension~$k$. But $\De0=\mbar{0,n}$ and $\De{d}$ (a collection of
points) are nonsingular. All $\Do{k}$ are nonsingular. Of course
$\Do{0}=\m{0,n}$ and $\Do{d}=\De{d}$. We have the long exact
sequence
$$\dots\to H^{a-1}_c(\De{k+1})\to H^a_c(\Do{k})\to H^a_c(\De{k})\to
H^a_c(\De{k+1})\to \dots$$ of mixed Hodge structures. Since the
$\De{k}$ are invariant for the natural action of $\Sn$, it is also
a sequence of $\Sn$-representations.

\begin{lm}
\label{rows0}
The cohomology groups $H^i(\m{0,n})$ vanish for $i>n-3$. For $0\le i\le n-3$,
the irreducible representations of $\Sn$ occurring in $H^i(\m{0,n})$
have Young diagrams with at most $i+1$ rows.
In particular, the irreducible representations of $\Sn$ occurring in
$H^\bullet(\m{0,n})$ have Young diagrams with at most $n-2$ rows.
\end{lm}

\bpf The claimed vanishing is immediate. Let us abbreviate the
rest of the statement by ``$H^i(\m{0,n})$ has $\le i+1$ rows''. We
prove it by induction on $n$. The case $n=3$ is trivial. Assume
$n>3$. 
Recall that $d=n-3$. 
We require an analysis of the boundary strata:

\smallskip \noindent
{\bf Claim.} Assume $d-b>0$. Then $H^a_c(\De{d-b})$ has $\le
d+1+b-a$ rows.

\smallskip
\noindent We prove the claim by induction on $b$. We begin with
the case $b=0$. Since $\De{d}$ is a collection of points, $a=0$ may
be assumed. Each point corresponds to a stable curve with $d$
nodes, hence with $d+1$ components. Each component has exactly
three special points (nodes or marked points). Let $n_j$ be the
number of marked points on the $j$th component, for some numbering
of the components. By permuting the $n$ marked points on the
stable curve, we obtain a $\Sn$-representa\-tion $R$, which is a
direct summand of $H^0(\De{d})$. Note that $R$ is a
subrepresentation of the induced representation
$$\Ind^{\Sn}_{\prod_{j=1}^{d+1}\Sigma_{n_j}}{\bf 1}.$$
The induced representation has $\le d+1$ rows, hence $R$ does. Now
$H^0(\De{d})$ is a direct sum of representations analogous to $R$,
thus it has $\le d+1$ rows as well. This proves the claim in the
case $b=0$.

Assume $b>0$. Observe that $H^a_c(\Do{k})\cong
H^{2(d-k)-a}(\Do{k})$ as $\Sn$-representa\-tions. Also, each
connected component of $\Do{k}$ is for $k\ge1$ a product of $k+1$
spaces $\m{0,m_j}$, with $m_j<n$. By induction on $n$ and the
K\"unneth formula, $H^a_c(\Do{k})$ has $\le 2(d-k)-a+k+1=2d-k-a+1$
rows, for $k\ge1$. Putting $k=d-b$, we find that $H^a_c(\Do{d-b})$
has $\le d+1+b-a$ rows.

By induction on $b$, we have that $H^a_c(\De{d-b+1})$ has $\le
d+b-a$ rows. From the long exact sequence, we find that
$H^a_c(\De{d-b})$ has $\le d+1+b-a$ rows. This proves the claim.

In particular, $H^a_c(\De1)$ has $\le 2d-a$ rows. Consider the
exact sequence
$$ H^{k-1}_c(\De1)\to H^k_c(\m{0,n})\stackrel{\alpha}{\to}
H^k_c(\mbar{0,n}).$$ From Lemma~\ref{wt2i} we know that
$H^k_c(\m{0,n})$ has weight $2k-2d$. But $H^k_c(\mbar{0,n})$ has
weight $k$. Thus $\alpha=0$ for $k<2d$. Hence $H^k_c(\m{0,n})$ has
$\le 2d+1-k$ rows for $k<2d$. Thus $H^i(\m{0,n})$ has $\le i+1$
rows for $i>0$. But it is obviously true for $i=0$ as well. This
finishes the proof.\epf

\section{The contribution of the boundary}\label{bdy}
In this section we determine the contribution of the boundary
$$\bm{1,n}=\mbar{1,n}\setminus\m{1,n},$$
i.e., we determine
$$\lan s_{1^n},e_c^{\Sn}(\bm{1,n})\ran.$$
We use the main result of \cite{Ge2}. 
To state it, we introduce the following notations:
\begin{eqnarray*}
&&{\bf a}_g := \sum_{n>2-2g} \ch_n(e_c^{\Sn}(\m{g,n})),\quad{\text {and}}\\
&&{\bf b}_g := \sum_{n>2-2g} \ch_n(e_c^{\Sn}(\mbar{g,n})).
\end{eqnarray*}
Here $\ch_n$ 
denotes 
the characteristic of a finite-dimensional
$\Sn$-represen\-tation (\cite{GK},~7.1) and its extension by linearity
to virtual representations.
For a (formal) symmetric function $f$ (such as ${\bf a}_g$ and ${\bf b}_g$),
we also write
\begin{eqnarray*}
&&f'=\frac{\partial f}{\partial p_1}=p_1^\perp f,\\
&&\dot{f}=\frac{\partial f}{\partial p_2}=\frac12\,p_2^\perp f,\\
&&\psi_i(f)=p_i\circ f.
\end{eqnarray*}
Here $p_i$ is the symmetric function equal to the sum of the $i$th powers
of the variables, 
$p_i^\perp$ is the adjoint of multiplication with $p_i$ 
with respect to the standard inner product,
and $\circ$ is the plethysm of symmetric functions (\cite{GK},~7.2).
We will denote the $i$th complete symmetric function by~$h_i$ and
the $i$th elementary symmetric function by~$e_i$.

We can now state Getzler's result (Theorem 2.5 in \cite{Ge2}):
$$
{\bf b}_1= \biggl( {\bf a}_1 - \frac{1}{2} \sum_{n=1}^\infty
\frac{\phi(n)}{n} \log(1-\psi_n({\bf a}_0'')) +
\frac{{\dot{{\bf a}}_0}^2+\dot{{\bf a}}_0+\tfrac14\psi_2({\bf a}_0'')}
{1-\psi_2({\bf a}_0'')} \biggr) \circ
(h_1+{\bf b}_0').
$$
The numerator of the third term inside the big
parentheses on the right-hand side
has been corrected here; there is a minor computational mistake in the
derivation of the theorem in line 4 on page 487, which affects the result
(but not Corollary 2.8).    

As Getzler remarks, the term $\bA1\circ(h_1+\bB0')$
corresponds to the sum over graphs obtained by attaching a forest
whose vertices have genus~$0$ to a vertex of genus~$1$; in particular,
$\bA1\circ h_1=\bA1$, the contribution of smooth curves, is part of this term,
corresponding to graphs consisting of a single vertex of genus~$1$.
The remainder of this term, corresponding to graphs where at least one vertex
of genus~$0$ has been attached to a vertex of genus~$1$, is part of the
contribution of the boundary. We show that the alternating representation
does not occur here. For a symmetric function
$$f=\sum_{n=0}^\infty f_n,$$
we write
$$\Alt(f)=\sum_{n=0}^\infty \lan s_{1^n},f_n\ran t^n.$$

\begin{lm}
\label{notrees1}
The alternating representation does not occur in the 
contribution of the part of the boundary
of~$\m{1,n}$ corresponding to graphs where at least one vertex of
genus~$0$ has been attached to a vertex of genus~$1$. In terms of the
notation introduced above:
$$\Alt(\bA1\circ(h_1+\bB0'))=\Alt(\bA1).$$
\end{lm}

\bpf
We choose to give a somewhat geometric proof instead of a proof using mostly
the language of symmetric functions.

Observe first that a boundary stratum corresponding to a graph with a
genus~$1$ vertex is isomorphic to a product
$$\m{1,m}\times\prod_i \m{0,n_i};$$
i.e., it is not necessary to take the quotient by a finite group.
(The corresponding graph has no automorphisms: there is a unique
shortest
path from each of the $n$ legs to the vertex of genus~$1$, and
every vertex and every edge lie on such a path.)
By the K\"unneth formula,
the cohomology of such a product is isomorphic to the tensor product
of the cohomologies of the factors.

Consider the $\Sn$-orbit of such a stratum. The direct sum of the cohomologies
of the strata in the orbit forms a $\Sn$-representation $V$.
It is induced from the cohomology of a single stratum, considered as a
representation $W$ of the stabilizer $G$ in $\Sn$ of the stratum.
By Frobenius Reciprocity, $V$ contains a copy of the
alternating representation if and only if $W$ contains a copy of the
restriction of the alternating representation to $G$.

To each vertex of the graph, one associates the symmetric group corresponding
to the legs attached to the vertex.
The product over the vertices of these symmetric groups is a subgroup $H$
of $G$ and the further restriction of the alternating representation to
$H$ is the tensor product over the vertices of the alternating representations
of these symmetric groups.

Consider a moduli space $\m{0,k}$ corresponding to an extremal vertex of
the graph corresponding to a boundary stratum as above.
The symmetric group associated to this vertex is a standard
subgroup $\Sigma_{k-1}\subset\Sigma_k$, permuting the $k-1$ legs
attached to the vertex and leaving the unique half-edge fixed.
By Lemma~\ref{rows0}, the irreducible $\Sigma_k$-representations occurring
in $H^\bullet(\m{0,k})$ have Young diagrams with at most $k-2$ rows.
The Young diagrams of the irreducible representations occurring in the
restriction to $\Sigma_{k-1}$ also have at most $k-2$ rows, as they are
obtained by removing one box. Therefore the alternating representation
does not occur here. It follows that $V$ does not contain a copy of the
alternating representation either.
\epf

We return to Getzler's result. We need to evaluate
$$\Alt \biggl( \biggl( - \frac{1}{2} \sum_{n=1}^\infty
\frac{\phi(n)}{n} \log(1-\psi_n({\bf a}_0'')) +
\frac{{\dot{{\bf a}}_0}^2+\dot{{\bf a}}_0+\tfrac14\psi_2({\bf a}_0'')}
{1-\psi_2({\bf a}_0'')} \biggr) \circ (h_1+{\bf b}_0') \biggr).
$$
Getzler remarks that the two terms inside the big inner parentheses may
be thought of as a sum over necklaces (graphs consisting of a single
circuit) and a correction term, taking into account the fact that
necklaces of~$1$ or~$2$ vertices have non-trivial involutions (while
those with more vertices do not).
The plethysm with $h_1+\bB0'$ stands again for attaching a forest whose
vertices have genus~$0$. We begin with the analogue of Lemma~\ref{notrees1}.

\begin{lm}
\label{notrees0}
The alternating representation does not occur in the
contribution of the part of the boundary
of~$\m{1,n}$ corresponding to graphs where at least one vertex of
genus~$0$ has been attached to a necklace. In terms of the
notation introduced above:
$$
\Alt \biggl( \biggl( - \frac{1}{2} \sum_{n=1}^\infty
\frac{\phi(n)}{n} \log(1-\psi_n({\bf a}_0'')) +
\frac{{\dot{{\bf a}}_0}^2+\dot{{\bf a}}_0+\tfrac14\psi_2({\bf a}_0'')}
{1-\psi_2({\bf a}_0'')} \biggr) \circ (h_1+{\bf b}_0') \biggr)
$$$$
=\Alt \biggl( - \frac{1}{2} \sum_{n=1}^\infty
\frac{\phi(n)}{n} \log(1-\psi_n({\bf a}_0'')) +
\frac{{\dot{{\bf a}}_0}^2+\dot{{\bf a}}_0+\tfrac14\psi_2({\bf a}_0'')}
{1-\psi_2({\bf a}_0'')} \biggr).
$$
\end{lm}

\bpf
In this case, each boundary stratum is isomorphic to a product
$$\biggl(\,\prod_{v\in\text{necklace}}\m{0,n(v)}\biggr)\big/I
\times \prod_{v\notin\text{necklace}}\m{0,n(v)}.
$$
The finite group $I$ is trivial when the necklace has at least $3$ vertices.
It has $2$ elements when the necklace has $1$ resp.~$2$ vertices and acts
by reversing the edge in the necklace resp.~by interchanging the two
edges of the necklace. In particular, $I$ acts trivially on the moduli
spaces corresponding to the vertices of the forest.

Just as in the proof of Lemma~\ref{notrees1}, the alternating representation
does not occur in the cohomology of a moduli space $\m{0,k}$
corresponding to an extremal vertex of one of the trees of the forest.
It follows that the alternating representation does not occur in the
cohomology of a $\Sn$-orbit of boundary strata as soon as the forest
is nonempty.
\epf

In order to determine the contribution of the part of the boundary
of~$\m{1,n}$ corresponding to necklaces without attached trees, we need
several lemmas.

\begin{lm}
\label{a0''}
The restriction of the $\Sn$-representation $H^\bullet(\m{0,n})$
to the standard subgroup~$\Sigma_{n-2}$ contains the alternating
representation exactly once. In terms of the notation introduced above: 
$$\Alt(\bA0'')=\frac{t}{1+t}\,.$$
\end{lm}

\bpf
The Young diagrams corresponding to the irreducible representations
of $\Sigma_{n-2}$ occurring in $\bA0''$ are
obtained by removing $2$ boxes from a Young diagram occurring in
$\bA0$. To obtain a copy of the alternating representation of
$\Sigma_{n-2}$, one needs to start with a Young diagram with at least
$n-2$ rows. From Lemma~\ref{rows0}, only the top cohomology
$H^{n-3}(\m{0,n})$ can contribute.
Observe that $H^{n-3}(\m{0,n})\cong H_c^{n-3}(\m{0,n})\cong H_{n-3}(\m{0,n})$
as $\Sigma_n$-representations.
Thus
\begin{eqnarray*}
\Alt(\bA0'')&=&\Alt\biggl(\frac{\p^2}{\p p_1^2}\sum_{n=3}^\infty
\ch_n\bigl(e_c^{\Sn}(\m{0,n})\bigr)\biggr) \\
&=&\Alt\biggl(\frac{\p^2}{\p p_1^2}\sum_{n=3}^\infty
(-1)^{n-3}\ch_n\bigl(H_c^{n-3}(\m{0,n})\bigr)\biggr).
\end{eqnarray*}
Getzler shows in \cite{Ge1}, p.~213, l.~3 that
$$H_c^{n-3}(\m{0,n}) \cong \sgn_n\otimes \Lie\(n\).$$
Here $\sgn_n$ 
denotes 
the alternating representation and $\Lie\(n\)$ 
the
$\Sn$-represen\-tation that is part of the cyclic Lie operad.
Getzler and Kapranov show in \cite{GK}, Example 7.24 that
$$
\Ch(\Lie):=\sum_{n=3}^\infty\ch_n(\Lie\(n\))
=(1-p_1)\sum_{n=1}^\infty\frac{\mu(n)}{n}\log(1-p_n)+h_1-h_2,
$$
where $\mu(n)$ is the M\"obius function.
Hence
$$
\sum_{n=3}^\infty(-1)^{n-3}\ch_n(\sgn_n\otimes\Lie\(n\))
=-(1+p_1)\sum_{n=1}^\infty\frac{\mu(n)}{n}\log(1+p_n)+h_1+e_2$$
and
\begin{eqnarray*}
&&
\frac{\p^2}{\p p_1^2}
\biggl(\,\sum_{n=3}^\infty(-1)^{n-3}\ch_n(\sgn_n\otimes\Lie\(n\))\biggr)
\\&&\qquad
=\frac{\p}{\p p_1}
\biggl(-\sum_{n=1}^\infty\frac{\mu(n)}{n}\log(1+p_n)-
(1+p_1)\frac{1}{1+p_1}\biggr)+1
\\&&\qquad
=1-\frac{1}{1+p_1}
=\frac{p_1}{1+p_1}\,.
\end{eqnarray*}
But
$$\Alt\biggl(\frac{p_1}{1+p_1}\biggr)=\frac{t}{1+t}\,,$$
since $\lan p_1^n,s_{1^n}\ran=1$.
\epf

\begin{lm}
\label{pk}
Let $f_n$ be a symmetric function of degree $n$. Assume that
$\lan s_{1^n},f_n\ran=0$. Then $\lan s_{1^{nk}},p_k\circ f_n\ran=0$.
\end{lm}

\bpf
Write $e(\lambda)$ resp.~$o(\lambda)$ for the number of even resp.~odd
parts of a partition $\lambda$.
For $\lambda$ a partition of $n$,
$$ o(\lambda)\equiv n\pmod{2} \qquad{\rm and}\qquad
\lan s_{1^n},p_{\lambda}\ran=(-1)^{e(\lambda)}.$$
Here $p_{\lambda}=\prod_i p_{\lambda_i}$ is the symmetric function
of degree~$n$ that is the product of the power sums corresponding
to the parts of~$\lambda$.
If $f_n=\sum_{\lambda} a_{\lambda}p_{\lambda}\,$, then $p_k\circ f_n
=\sum_{\lambda} a_{\lambda}p_{k\lambda}\,$, where $k\lambda$ is
the partition of $kn$ obtained from $\lambda$ by multiplying all parts
with $k$. For $k$ odd,
$$\sum_{\lambda} a_{\lambda}(-1)^{e(\lambda)}
=\sum_{\lambda} a_{\lambda}(-1)^{e(k\lambda)},$$
whereas for $k$ even,
$$e(k\lambda)=o(\lambda)+e(\lambda)\equiv n+e(\lambda) \pmod{2},$$
so that
$$\sum_{\lambda} a_{\lambda}(-1)^{e(k\lambda)}=(-1)^n
\sum_{\lambda} a_{\lambda}(-1)^{e(\lambda)}.$$
The result follows.
\epf

\begin{lm}
\label{psik}
The occurrence of the alternating representation observed in
Lemma~\ref{a0''} is stable under plethysm with~$p_k$. In terms
of the notation introduced above:
$$\Alt(\psi_k(\bA0''))=\frac{-(-t)^k}{1-(-t)^k}\,.$$
\end{lm}

\bpf
{From} Lemma~\ref{a0''},
$$\Alt(\bA0'')=\Alt\biggl(\frac{p_1}{1+p_1}\biggr)\,.$$
Applying Lemma~\ref{pk} and using
that $\lan s_{1^{kn}},p_k^n\ran=(-1)^{(k-1)n}$, we find
$$
\Alt(\psi_k(\bA0''))=\Alt\biggl(p_k\circ\biggl(\frac{p_1}{1+p_1}\biggr)
\biggr)=\Alt\biggl(\frac{p_k}{1+p_k}\biggr)
$$$$ 
=\Alt\biggl(\sum_{n=1}^\infty(-1)^{n-1}p_k^n\biggr)
=\sum_{n=1}^\infty(-1)^{n-1}(-1)^{(k-1)n}t^{kn}
$$$$ 
=-\sum_{n=1}^\infty ((-t)^k)^n
=\frac{-(-t)^k}{1-(-t)^k}\,.
$$
\epf

\noindent For two symmetric functions $f$ and $g$, we have
$$\Alt(fg)=\Alt(f)\Alt(g).$$
This follows immediately from the Littlewood-Richardson rule
(cf.~\cite{FH}, p.~456). One may use this identity to shorten the proof
of Lemma~\ref{psik}. Similarly,
$$
\Alt(\log(1-\psi_k(\bA0'')))=\log\biggl(1-\frac{-(-t)^k}{1-(-t)^k}\biggr)
=-\log(1-(-t)^k)$$
and
$$
\Alt\biggl(-\sum_{n=1}^{\infty}\frac{\phi(n)}{n}\log(1-\psi_n(\bA0'')\biggr)
=\sum_{n=1}^{\infty}\frac{\phi(n)}{n}\log(1-(-t)^n)
$$$$
=-\sum_{n=1}^{\infty}\frac{\phi(n)}{n}\sum_{k=1}^{\infty}\frac{(-t)^{nk}}{k}
=-\sum_{m=1}^{\infty}\frac{(-t)^m}{m}\sum_{d|m}\phi(d)
$$$$ 
=-\sum_{m=1}^{\infty}(-t)^m=\frac{t}{1+t}\,.
$$
It remains to evaluate the contribution of the correction term,
$$
\Alt\biggl(\frac{{\dot{{\bf a}}_0}^2+\dot{{\bf a}}_0
+\tfrac14\psi_2({\bf a}_0'')}{1-\psi_2({\bf a}_0'')} \biggr).
$$
We need one more lemma.

\begin{lm}
\label{a0dot}
The alternating part of the formal symmetric function
$$\dot{{\bf a}}_0=\frac{\p {\bf a}_0}{\p p_2}$$ is given
by the following formula:
$$\Alt(\dot{{\bf a}}_0)=\frac12\,\frac{t}{1-t}\,.$$
\end{lm}

\bpf
In terms of Young diagrams, multiplication by $s_2$ is the operation
of adding two boxes, not in the same column, and multiplication by $s_{1^2}$
is the operation of adding two boxes, not in the same row. Now
$p_2=s_2-s_{1^2}$ and $\tfrac{\p}{\p p_2}=\tfrac12p_2^{\perp}$, where
$p_2^{\perp}$ is the adjoint of multiplication by $p_2$. Thus, to obtain
a copy of the alternating representation of $\Sigma_{n-2}$ in a term of
$\dot{{\bf a}}_0$, one needs to start with a Young diagram with at least
$n-2$ rows, just as in the proof of Lemma~\ref{a0''}. So
$$\Alt(\dot{{\bf a}}_0)=\Alt\biggl(\frac{\p}{\p p_2}\sum_{n=3}^\infty
(-1)^{n-3}\ch_n\bigl(H_c^{n-3}(\m{0,n})\bigr)\biggr).$$
We now find
$$
\frac{\p}{\p p_2}
\biggl(\,\sum_{n=3}^\infty(-1)^{n-3}\ch_n(\sgn_n\otimes\Lie\(n\))\biggr)
=\frac12\,\frac{1+p_1}{1+p_2}-\frac12
=\frac12\,\frac{p_1-p_2}{1+p_2}\,.
$$
But
$$\Alt\biggl(\frac12\,\frac{p_1-p_2}{1+p_2}\biggr)
=\frac12\,\frac{t+t^2}{1-t^2}=\frac12\,\frac{t}{1-t}\,.$$
\epf

\noindent An easy calculation combining Lemmas~\ref{psik}
and~\ref{a0dot} gives
$$
\Alt\biggl(\frac{{\dot{{\bf a}}_0}^2+\dot{{\bf a}}_0
+\tfrac14\psi_2({\bf a}_0'')}{1-\psi_2({\bf a}_0'')} \biggr)
=\frac12\,\frac{t}{1-t}\,.
$$
The contribution from the necklaces becomes then
$$
\frac12\,\frac{t}{1+t}+\frac12\,\frac{t}{1-t}=\frac{t}{1-t^2}\,,$$
i.e., $1$ for $n$ odd and $0$ for $n$ even.
We have proved the following result.

\begin{tm}
\label{boundary}
The alternating part of the $\Sn$-equivariant Euler characteristic
of the cohomology of~$\bm{1,n}$ is given by the
following formula:
$$\lan s_{1^n},e_c^{\Sn}(\bm{1,n})\ran=\left\{
\begin{array}{rl}
1,&\quad n\,\,{\rm odd;}\\
0,&\quad n\,\,{\rm even.}\\
\end{array}
\right.
$$
\end{tm}

\section{The construction of the motive}
The main result of Section~\ref{bdy} (Theorem~\ref{boundary}) is
$$A_c(\bm{1,n})=\left\{
\begin{array}{rl}
1,&\quad n\,\,{\rm odd;}\\
0,&\quad n\,\,{\rm even.}\\
\end{array}
\right.
$$
Combining this with Theorem~\ref{interior}, we immediately obtain
$$A_c(\mbar{1,n})=\left\{
\begin{array}{rl}
-S[n+1],&\quad n\,\,{\rm odd;}\\
0,&\quad n\,\,{\rm even.}\\
\end{array}
\right.
$$
Let $n>1$ be an odd integer. The pair consisting of
$\mbar{1,n}$ and the projector
$$\frac1{n!}\sum_{\sigma\in\Sn}(-1)^{\sgn(\sigma)}\sigma_*$$
defines a Chow motive, since $\mbar{1,n}$ is the quotient of a 
smooth projective variety by a finite group~\cite{BP}.
We wish to show that it is pure of
degree $n$. This will conclude our construction of the motive
$S[n+1]$, an alternative (in level $1$ only) to Scholl's
construction. The arguments below are similar to those
in~\cite{Sc}, 1.3.4.

First, $\lan s_{1^n},H^i_c(\m{1,n})\ran=0$ for $i\neq n$. As in~\cite{De},
proof of~5.3, the degeneration of the Leray spectral sequence at $E_2$ due
to Lieberman's trick implies that $\lan s_{1^n},H^i_c(\E^{n-1})\ran=0$
when $i\neq n$ for a relative elliptic curve $\E\to S$ and this implies
the statement.

Thus we have an exact sequence
$$0 \to H_c^{n-1}(\mbar{1,n})(\alpha) \to H_c^{n-1}(\bm{1,n})(\alpha)
\to H^n_c(\m{1,n})(\alpha) \to $$
$$\to H^n_c(\mbar{1,n})(\alpha) \to H^n_c(\bm{1,n})(\alpha) \to 0$$
and isomorphisms $H^i(\mbar{1,n})(\alpha)\cong
H^i(\bm{1,n})(\alpha)$ for $i\notin\{n-1,n\}$. Here $V(\alpha)$
denotes the alternating part of a $\Sn$-representation $V$.
Therefore $H^i(\bm{1,n})(\alpha)$ is pure of weight $i$ for $i>n$.
But then all these spaces vanish, since $H^i(\bm{1,n})$ has weight
$\le i$ for all $i$ and since $A_c(\bm{1,n})$ has weight $0$.
Hence $H^i(\mbar{1,n})(\alpha)=0$ for $i>n$ and then by duality
for $i<n$ as well. This shows that
$H^n(\mbar{1,n})(\alpha)=S[n+1]$ and concludes the alternative
construction of these motives.

\end{document}